\documentclass[11pt,a4paper]{article}

\usepackage[latin1]{inputenc}

\usepackage[normalem]{ulem}

\usepackage{amssymb}

\usepackage{graphicx}

\usepackage{graphicx,indentfirst,amsmath,amsfonts,amssymb,amsthm,newlfont}
\usepackage{epsfig}
\newtheorem{theorem}{Theorem}[section]
\newtheorem{proposition}[theorem]{Proposition}
\newtheorem{remark}[theorem]{Remark}
\newtheorem{definition}[theorem]{Definition}

\newtheorem{example}[theorem]{Example}

\newcommand{\dif}{\mathrm{Diff}}
\newcommand{\difH}{\mathrm{Diff}(\Delta^H, M)}
\newcommand{\difV}{\mathrm{Diff}(\Delta^V, M)}

\newcommand{\Prob}{\mathbf P}
\newcommand{\R}{\mathbf R}
\newcommand{\N}{\mathbf N}
\newcommand{\Z}{\mathbf Z}

\begin{document}
\begin{center}


 {\Large {\bf Extension of time for decomposition of stochastic flows in spaces with complementary foliations
 }}

\end{center}

\vspace{0.3cm}

\begin{center}
{\large { Leandro Morgado}\footnote{E-mail:
morgado@ime.unicamp.br.
Research supported by  FAPESP 11/14797-2.} \ \ \ \ \ \ \ \ \ \ \ \ \
{ Paulo R. Ruffino}\footnote{Corresponding author, e-mail:
ruffino@ime.unicamp.br.
Research partially supported by FAPESP nr. 12/18780-0, 11/50151-0 and CNPq 477861/2013-0.}}

\vspace{0.2cm}

\textit{Departamento de Matem\'{a}tica, Universidade Estadual de Campinas, \\
13.083-859- Campinas - SP, Brazil.}

\end{center}

\begin{abstract}
Let $M$ be a manifold equipped (locally) with a pair of complementary foliations. In Catuogno, da Silva and Ruffino \cite{CSR2}, it is shown that, up to a stopping time $\tau$, a stochastic flow of local diffeomorphisms $\varphi_t$ in $M$ can be decomposed in diffeomorphisms that preserves this foliations. In this article we present techniques which allows us to extend the time of this decomposition. For this extension, we use two techniques: In the first one, assuming that the vector fields of the system commute with each other, we apply Marcus equation to jump nondecomposable diffeomorphisms. The second approach deals with the general case: we introduce a `stop and go' technique that allows us to construct a process that follows the original flow in the `good zones' for the decomposition, and remains paused in `bad zones'. Among other applications, our results open the possibility of studying the asymptotic behaviour of each component.
\end{abstract}


\vspace{0.3cm}
\noindent {\bf MSC2010 subject classification:} 60H10, 60D05, (57R30).

\section{Introduction}

Given a stochastic flow $\varphi_t$ of local diffeomorphisms in a differentiable manifold $M$, in many circunstances
the decomposition of $\varphi_t$ with components in subgroups of
the group of diffeomorphisms $\dif(M)$ provide interesting dynamical or geometrical information of the
system.  In the literature, this kind of
decomposition has been studied in several  frameworks and with different
aimed subgroups; among others in Bismut \cite{Bismut}, Kunita
\cite{Kunita-1}, \cite{Kunita-2},  Ming Liao \cite{ML} and some of our previous
work \cite{CSR}, \cite{Colonius-Ruffino}, \cite{Ruffino}. In the
last few papers mentioned, geometrical conditions on a Riemannian manifold have
been stated to guarantee the existence of the decomposition where the first
component lies in the subgroups of
isometries or affine
transformations.

In particular, in Catuogno, da Silva and Ruffino \cite{CSR2}, the authors consider a pair of complementary distributions in a differentiable manifold $M$, in the sense that each tangent space splits in a direct sum of two subspaces depending differentiably on the points in $M$. These subspaces are called, by convenience, horizontal and vertical distributions. In that article it is shown that locally, up to a stopping time $\tau$, a stochastic flow $\varphi_t$ in $M$ can be decomposed as $\varphi_t=\xi_t \circ \Psi_t$, where $\xi_t$ is a diffusion in the group of diffeomorphisms $\difH$ generated by horizontal vector fields. On the other hand, $\Psi_t$ is a process in the group of diffeomorphisms $\difV $ generated by vertical vector fields. The infinite dimensional Lie group structure considered in this case is described in Milnor \cite{Milnor}, Neeb \cite{Neeb} and Omori \cite{Omori}.  In \cite{CSR2} the authors present 
stochastic differential equations on the corresponding infinite dimensional Lie subgroups for the components $\xi_t$ and $\Psi_t$.

The stopping time $\tau$ mentioned above, which restricts the time where the decomposition holds, appears due to an
explosion in the equation of one of the components of the decomposition, with initial condition at the identity map.
This explosion phenomenon is related to the choice of the distributions, and  it is not intrinsic to the original flow $\varphi_t$. Moreover, there might be large intervals on time after $\tau$ such that the diffeomosphisms $\varphi_t$
are decomposable for $t$ in these intervals.

In this article we present techniques which allows us to extend the time of the decomposition in a simplified framework of \cite{CSR2}: We consider that the distributions are in fact integrable, hence the manifold $M$ has locally a pair of complementary foliations. This is a natural structure since any coordinate systems on $M$ generates this foliations. Besides, it means that locally the manifold is diffeomorphic to an open set in $\R^n = \R^{p} \times \R^{n-p}$, where $p$ and $n-p$ are the dimensions of the horizontal and vertical foliations, respectively.

With this change of coordinates in mind, consider a diffeomorphism $\varphi: U\subset \R^{p} \times \R^{n-p} \rightarrow \varphi(U) \subset \R^{p} \times \R^{n-p}$, with $U$ an open set of $\R^{n}$ . One writes in coordinates $\varphi = (\varphi^1 (x,y), \varphi^2 (x,y))$. Then, by the inverse function theorem, the local decomposition exists if and only if the $(n-p)\times (n-p)$ matrix $\frac{\partial \varphi^2(x,y)}{\partial y}$ is invertible. That is, $ \varphi = (\xi^1 (x,y), Id_2) \circ (Id_1, \varphi^2 (x,y))$, where $Id_1$ and $Id_2$ are the identities in $\R^p$ and in $\R^{n-p}$, respectively. With this notation, we call

\begin{definition}
A diffeomorphism $\varphi: U\subset \R^{p} \times \R^{n-p} \rightarrow \varphi(U) \subset \R^{p} \times \R^{n-p}$ is $p$-decomposable in a neighbourhood of $(x,y)$ if $ \det \frac{\partial \varphi^2(x,y)}{\partial y}\neq 0$.
\end{definition}

The explosion time $\tau$ for the decomposition of the flow is related to the fact that $\varphi_{\tau}$ is not $p$-decomposable. As a basic example to illustrate the explosion, consider the linear rotation, whose $1$-decomposition is given by
\[
\varphi_t = \left(
                                \begin{array}{cc}
                                  \cos t & -\sin t \\
                                  \sin t & \cos t
                                \end{array}
                              \right)=
                              \left(
                                \begin{array}{cc}
                                  \sec t & -\tan t \\
                                  0 & 1
                                \end{array}
                              \right)
                              \left(
                                \begin{array}{cc}
                                  1 & 0 \\
                                  \sin t & \cos t
                                \end{array}
                              \right).
\]
Hence, $\varphi_t$ is not decomposable when $t=\frac{\pi}{2} + k \pi$, with $k\in \Z$. The stopping time as in \cite{CSR2} is given by $\tau = \frac{\pi}{2}$, nevertheless, the decomposition can be done in many other intervals after this point. Note that, in this example, immediately after this stopping time, $\cos t$ changes its sign and the second component $\Psi_t$ revert orientation, hence, it is no longer in the subgroup of diffeomorphism generated by vertical vector fields, but it is rather in the subgroup of diffeomorphisms which preserves each vertical leaves of the foliation, not necessarily in the connected component of the identity. Hence, here $\difH$ and $\difV$ are extended to include nonpreserving orientation diffeomorphisms.

One of the motivations for the decomposition in \cite{CSR2} is that it shows how close the system is horizontal leaf-preserving, i.e., how close is the vertical component to the identity. Here, the motivation is to show that this decomposition extends further in time, even if one has to consider the nonpreserving orientation diffeomorphisms. Among other applications, our results open the possibility of studying the asymptotic behaviour of each component.

 The extension of time is described using two techniques: the first one (Section 2) assumes that the vector fields of the system commute with each other. We apply Marcus equation to jump nondecomposable diffeomorphisms; essentially we construct a flow $\widetilde{\varphi}_t$ which is close to the original one and is p-decomposable for all $t \geq 0$. Approximation here holds except in a set of arbitrarily small probability (Theorem 2.3). The second approach (Section 3) deals with the general case: we introduce a `stop and go' technique that allows us to construct a process that follows the original flow in the `good zones' for the decomposition, and remains paused in `bad zones', i.e., close to nondecomposable diffeomorphisms (Proposition 3.2).


\section{Commuting vector fields}

\subsection{Marcus equation and preliminaries}


For reader's convenience, we recall the definition of Marcus equation (see e.g. Kurtz, Pardoux e Protter \cite{KPP}):

\begin{equation}\begin{array}{lll} \label{eq: define marcus}
dx_t = \displaystyle \sum_{i=0}^m X^i(x_t) \circ dZ^i_t.
\end{array}
\end{equation}
with initial condition $x(0) = x_0$. The integral form is given by
\begin{equation}\begin{array}{lll}
x_t  = x_0 + \displaystyle \sum_{i=0}^m \int_0^t X^i(X_s) \circ dZ^i_s =:x_0 + \displaystyle \int_0^t X(x_s) \circ dZ_s
\end{array}
\end{equation}
where $x_t$ is an adapted stochastic process taking values in $\R^d$, the integrator $\{Z^i_s : s \geq 0\}$  is a semimartingale with jumps and $X^i$ are vector fields in $\R^d$ for all $i \in \{0,1,\ldots,m\}$. The solution is interpreted as a stochastic process that satisfies the equation:

\begin{equation}\begin{array}{lll} \label{eq: interpreta marcus}
x_t = x_0 + \displaystyle \int_0^t X(x_{s_-}) \ dZ_s \ + \ \frac{1}{2} \displaystyle \int_0^t X'X(x_s) \ d[Z,Z]_s^c \ + \\
\hspace{1cm} + \displaystyle \sum_{0 < s \leq t} \left\{ \varphi (X \Delta Z_s, x_{s_-}) - x_{s_-} - X(x_{s_-}) \Delta Z_s \right\},
\end{array}
\end{equation}
in the following sense: the first term on the right hand side of equation $(\ref{eq: interpreta marcus} )$ is a standard Itô integral of the predictable process $\{X(x_{t_-})\}$ with respect to the semimartingale $Z_t$. The second term is a Stieltjes integral with respect to the continuous part of the quadratic variation of $Z_t$. In the third term:  $\varphi (X \Delta Z_s, x_{s_-})$ indicates the solution of the inicial value problem in $t = 1$ with respect to the vector field $X \Delta Z_s$, and initial condition $x_{s_-}$. Thus, the jumps of Marcus equation occurs in the direction of the deterministic flow given by the vector field $X \Delta Z_s$.

Regularity conditions on the vector fields $X$ implies that there exists a unique stochastic flow of diffeomorphisms $\varphi_t$ that is solution of equation (\ref{eq: define marcus}). Conditions on the derivatives of the vector fields guarantee that the flow exists for all $t\geq 0$. Moreover, for an embedded submanifold $M$ in an Euclidean space, if the vector fields of the equation (\ref{eq: define marcus}) are in $TM$ then, for each initial condition on $M$, the solution stays in $M$ a.s., see \cite{KPP}.

\medskip


We present a simple example that illustrates the technique in the proof of the next theorem.
The idea is to use Marcus equation in order to construct a new flow  $\widetilde{\varphi}_t$ , which is decomposable for all $t\geq 0$ and simultaneouly is close to the original flow  $\varphi_t$. Consider again the simple example of pure rotation flow
        $$  \varphi_t = \left(
                                \begin{array}{cc}
                                  \cos t & -\sin t \\
                                  \sin t & \cos t
                                \end{array}
                              \right),$$
which is solution of the linear equation $dx_t= A x_t\ dt$, where $A$ is skew-symmetric. We have to construct an integrator $Z_t$ in  such a way that the solution flow $\widetilde{\varphi}_t$ of Marcus equation $dx_t = Ax_t \circ dZ_t$ jumps over the nondecomposable rotations, with the additional property that $\widetilde{\varphi}_t$ differs from ${\varphi}_t$ for $t$ in a set with Lebesgue measure arbitrarily small.


Fix $\varepsilon > 0$ and choose a sequence of points $p_n$ just before the critical time for the decomposition:
\[
p_n \in (\frac{\pi}{2} + n \pi - \frac{\varepsilon}{2^n}, \frac{\pi}{2} + n \pi)
\]
for all nonnegative integer $n$. These points indicate the beginning of what we call in the proof of the theorem a `red zone' for the corresponding critical point. Then, define

\[
Z_t  = \begin{cases}
 p_n \ \ \ \ \  \mbox{if} \ \ t \in [p_n, (2n+1) \pi - p_n); \\
t \ \ \ \ \ \ \ \mbox{otherwise}.
\end{cases}
\]
$Z_t$ is a semimartingale with jumps, and the solution of $dx_t = Ax_t \circ dZ_t$ is given by

\[
\widetilde{\varphi}_t = \left(
  \begin{array}{cc}
    \cos (Z_t) & -\sin (Z_t) \\
    \sin (Z_t) & \cos (Z_t) \\
  \end{array}\right).
 \]
We have that  $\widetilde{\varphi}_t$ is 1-decomposable for all $t \geq 0$ and, by construction, it is arbitrarily close to $\varphi_t$.

\subsection{Main results}

In this section, we consider the flow $\varphi_t$ generated by a Stratonovich  SDE:
\[
 dx_t =  X^0 (x_t)\ dt + \displaystyle \sum_{i=1}^{m} X^i(x_t) \circ dB^i_t
\]
where $(B^1_t, \ldots, B^m_t)\in \R^m$ is a Brownian motion in a filtered probability space $(\Omega, \mathcal{F}, \mathcal{F}_t, \Prob )$;  $X^i$ for $i \in \{0, 1,2,\ldots, m\}$ are smooth vector fields in $\R^n$. We assume that the corresponding (deterministic) flows $\phi^{X^i}_ t$ commute with each other.

Given $u= (t_0, t_1, \ldots , t_m) \in \R^{m+1}$, let $\phi ({u})$ be the composition of the deterministic flows, that is:
$$\phi ({u})  = \phi^{X_0}_{t_0} \circ \phi^{X_1}_{t_1} \circ \ldots \circ \phi^{X_m}_{t_m}.$$
We write $(x,y) \in \R^p \times \R^{n-p} = \R^n$ and $\phi ({u})(x,y)= (\phi({u})^1, \phi ({u})^2 )$. For a given initial condition, consider $P = \{u \in \R^{m+1}: \det \frac{\partial \phi(u)^2}{\partial y} = 0 \}$, the set of undecomposable generated diffeomorphisms.

\begin{theorem}
\label{browniano1}
Assume that $P$ has zero Lebesgue measure. Then, given $\varepsilon > 0$ and $a > 0$, there exists a semimartingale $Z_t\in \R ^{m+1}$  such that the solution $\widetilde{\varphi}_t$ of SDE equation $dx_t = \displaystyle \sum_{i=0}^{m} X^i(x_t) \circ d Z^i_t$ can be decomposed for $t \geq 0$. Besides, $C(\omega) = \{t \geq 0: \varphi_t(\omega) \neq \widetilde{\varphi}_t (\omega)\}$ is Lebesgue measurable a.s. and $\Prob [\mu(C) > a ] \leq \varepsilon$.

\begin{proof}

Inicially, fix $\varepsilon > 0$ and $a > 0$. For each $t \geq 0$, denote $U_t = (t, B_t^1, \ldots, B_t^n)$. As $P$ has zero Lebesgue measure, given $\delta > 0$, there is an open set $A_\delta$, such that $A_\delta \supseteq P$ and $\mu(A_\delta) < \delta$. By properties of Brownian motion, given $d > 0$, we have that:

$$\Prob \left[ \ U_s \in A_\delta \ \forall s \in [t, t + d] \ \big| \ U_s \in A_\delta \right] \to 0 \ \ \text{when} \ \ \delta \to 0.$$
Therefore, for each $k \in \N$, there exists an open set $A_k$ that contains $P$ with Lebesgue measure sufficiently small such that:
$$\Prob \left[ \ U_s \in A_k \ \forall s \in [t, t + \frac{a}{2^k}] \ \big| \ U_s \in A_k \ \right] < \frac{\varepsilon}{2^k}.$$

We can assume that the sequence of sets $\big( A_k \big)_{k \in \N}$ is decreasing;
besides, since $0 \notin P$ one can always assume that $0 \notin A_k$ for all $k \in \N$.

By Urysohn's lemma, there exist continuous functions $F_k: \R^m \to [0,1]$, such that $F_k^{-1} (1) = (A_k)^C$  and $F_k^{-1} (0) = P$. For each $k \in \N$, consider the following partition in $\R^{m+ 1}$:
\begin{itemize}
\item `Green zone' $ G_k := F_k^{-1} (1)$;
\item `Yellow zone' $ Y_k = F_k^{-1} \big(\frac{1}{2}, 1 \big)$;
\item `Red zone'  $\ R_k= F_k^{-1} \big[0 , \frac{1}{2}\big].$
\end{itemize}
Red zone corresponds to a set where we do not allow the dynamics to get into (via stopping times). Green zone corresponds to the set where we allow the dynamics to go on freely; finally, the yellow zone is an intermediate set.

\begin{figure}[h!]
\begin{center}
 \includegraphics[scale=0.6]{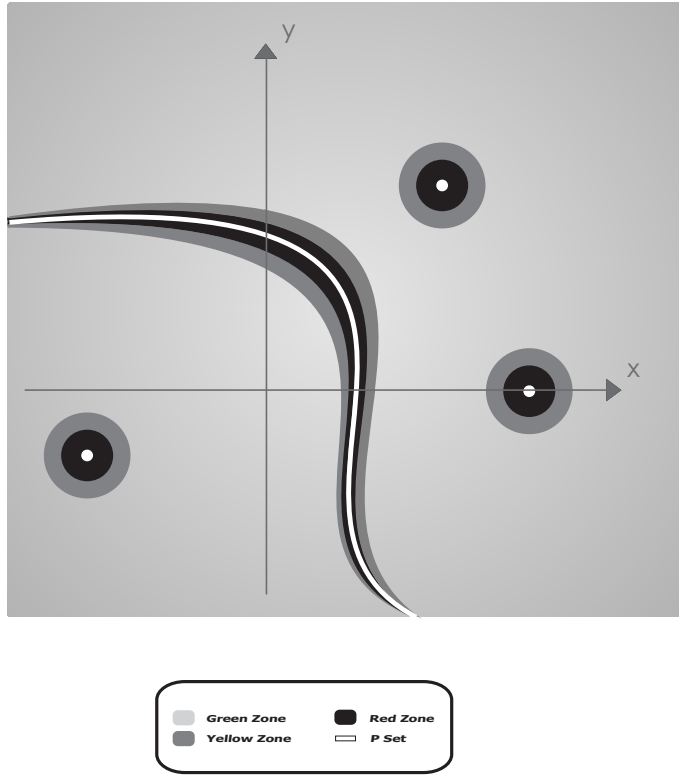}
 \caption{Sketch of green, yellow and red zones in $\R^2$.}
\end{center}
\end{figure}

By construction, the origin of $\R^{m+1}$ belongs to green zone for all $k \in \N$, besides, this set is increasing with $k$. We define by induction the following sequences of stopping times: Take $\overline{T}_0 = 0$, $T_i = \inf \{ t > \overline{T}_{i-1}: U_t \in R_i\}$ and $\overline{T}_i = \inf \{ t > T_i: U_t \in G_i\}$. The yellow zone guarantee that $ \overline{T}_{i-1} < T_i < \overline{T}_{i}$  for all $i \in \N$.

Define the semimartingale $\{Z_t: t \geq 0\}$ as follows:

\[
Z_t  = \begin{cases}
 U_t, \ \ \, \mbox{ if } \ \ t \in \big[\overline{T}_{k-1}, T_k \big) \ \text{ for some } \ k \in N; \\
U_{T_k}, \ \ \mbox{if } \ \ t \in \big[T_{k}, \overline{T}_k \big).
\end{cases}
\]
Note that $Z_t$ is constant when $U_t$ cross the yellow zone when going from the red to the green areas; otherwise, it coincides with  $U_t$. Now, consider the SDE given by
\[
dx_t = \displaystyle \sum_{i=0}^{m} X^i(x_t) \circ d Z^i_t,
\]
whose solution flow is given by:

\[
\widetilde{\varphi}_t = \phi_{Z_t} = \phi^{X_0}_{Z^0_t} \circ \phi^{X_1}_{Z^1_t} \circ \ldots \circ \phi^{X_m}_{Z^m_t}.
\]
Hence,  $\widetilde{\varphi}_t$ is $p$-decomposible for all $t \geq 0$.

For the second part of the statement, initially note that $C(\omega)$ is Lebesgue measurable a.s. since $\varphi_t$ and $\widetilde{\varphi}_t$ are measurable with respect to the product $\sigma$-algebra in $\R \otimes \Omega$. Moreover,   $\varphi_t (\omega) = \widetilde{\varphi}_t (\omega)$ if  $U_t (\omega) = Z_t (\omega)$, therefore:
$$C(\omega) \subseteq \displaystyle \bigcup_{k \in \N} \big(T_k(\omega), \overline{T}_k (\omega) \big).$$



\

Finally, a necessary condition to $\mu(C(\omega)) > a$ is that, for some $k \in \N$, we have $\overline{T}_k(\omega) - T_k (\omega) > \frac{a}{2^k}$. So, it follows that:

\begin{equation}\begin{array}{lll}
\Prob [ \mu(C) > a ] & \leq & \Prob \left[   \displaystyle \bigcup_{k \in \N} \left\{ \omega:  \big(\overline{T}_k - T_k \big) >  \frac{a}{2^k} \right\} \right] \vspace{3mm} \\
\vspace{3mm}
&\leq& \displaystyle \sum_{k \in \N} \Prob \left[ \big(\overline{T}_k - T_k \big) > \frac{a}{2^k} \right] \\
 &\leq& \displaystyle \sum_{k \in \N} \frac{\varepsilon}{2^k} = \varepsilon.
\end{array}
\end{equation}
\end{proof}

\end{theorem}

\begin{remark} {\em Generically in the $C^1$ topology, when the SDE has just a single vector field, then  the set $P\subset \R$ in the statement of the theorem has zero Lebesgue measure: In fact, the condition on the second  component $\det \frac{\partial \phi^2_t}{\partial y} = 0$ is destroyed by small perturbation on their derivatives.
Also for linear vector fields, $P \subset \R$ is discrete since $\det \frac{\partial \phi^2_t}{\partial y}$ is a nonvanishing analytical function.}
\end{remark}



When $P$, defined above, has an arbitrary Lebesgue measure, a weaker result holds. One can still construct a semimartingale $Z_t$ that is close to the process $U_t$, such that the corresponding generated flows $\widetilde{\varphi}_t$ and $\varphi_t$ are again close of each other. But here one loses the probabilistic approach to control the Lebesgue measure of $C(\omega)$. With the same notation as before:

%
%


 \medskip

\begin{proposition}
\label{fluxosquecomutam}
 Given $\varepsilon > 0$ and an open set $A \supseteq P$ such that $\mu (A \setminus P)< \varepsilon$, there exists a semimartingale $Z_t$ which satisfies the following properties: 
\begin{enumerate}
 \item The corresponding solution flow $\widetilde{\varphi}_t$ driven by $Z_t$
is $p$-decomposable for all $t \geq 0$;

\item If $U_t (\omega) \notin A$, then $Z_t(\omega)= U_t(\omega)$, hence $\widetilde{\varphi}_t (\omega)= \varphi_t (\omega)$.

\end{enumerate}

\begin{proof}


Since $0 \notin P$, we can assume (reducing $A$, if necessary), that $0 \notin A$.
By Urysohn's lemma, there exists a continuous function $F: \R^{m+1} \to [0,1]$, such that $F^{-1}(1) = A^C$ and $F^{-1}(0) = P$. In an analogous way to the proof of Theorem \ref{browniano1}, consider a partition of $\R^{m+1}$ as follows:
\begin{itemize}
\item `Green zone'  $ G = F^{-1} (1)$;
\item `Yellow zone' $Y = F^{-1} \big(\frac{1}{2}, 1 \big)$;
\item `Red zone'  $ R = F^{-1} \big[0 , \frac{1}{2}\big].$
\end{itemize}

We remark that here, different from the proof of the previous theorem, the green, yellow and red zones remain fixed during all the process. Again, by induction, consider the following sequences of stopping times: Take $\overline{T}_0 = 0$, $T_i = \inf \{ t > \overline{T}_{i-1}: U_t \in R\}$ and $\overline{T}_i = \inf \{ t > T_i: U_t \in G\}$. Define the process $\{Z_t: t \geq 0\}$ as follows:

\[
Z_t = \begin{cases}
U_t, \ \ \mbox{ if }  \ \ t \in \big[\overline{T}_{k-1}, T_k \big) \ \mbox{ for some } \ k \in N; \\
U_{T_k}, \ \ \mbox{if} \ \ t \in \big[T_{k}, \overline{T}_k \big).
\end{cases}
\]
Finally, note that if $U_t (\omega) \notin A$, then $U_t (\omega)$ is contained in the green zone. Therefore, in this case  $U_t (\omega) = Z_t (\omega)$ and $\varphi_t (\omega) = \widetilde{\varphi}_t (\omega)$.

\end{proof}
\end{proposition}

\begin{example} {\em
We present an example of $P$ with nonzero Lebesgue measure. Let $f: \R \rightarrow [-\pi/2, \pi/2]$ be a smooth function such that
$f(x)=-\pi/2$ if $x\leq - \pi/2$ and $f(x)=\pi/2$ if $x\geq \pi/2$. Consider the
following deterministic ODE in $\R^3$:
\[
\left( \begin{array}{c}
         x' \\
         y'\\
         z'
        \end{array}
\right) =
\left( \begin{array}{c}
         -f'(z)y \\
         f'(z)x\\
         1
        \end{array}
\right)
\]
whose solution at points in the horizontal plane $(x_0, y_0, 0)$ is given by

\[
\varphi(t)(x_0, y_0, 0)  =  \left(
\begin{array}{cc}
\cos f(t) & \begin{array}{cc} -\sin f(t) & 0  \end{array} \\
\begin{array}{c}  \sin f(t) \\ 0  \end{array} &
\put(-40,-10){\dashbox{2}(75,30){$ \begin{array}{cc} \ \ \ \cos f(t) & \  0 \\ 0 & \   1 \end{array} $}}
\end{array} \right) \left( \begin{array}{c}
                           x_0 \\
                             y_0\\
                             t
                            \end{array} \right)  .
\]

Consider the $1$-decomposition of the flow in  $\R\times \R^2 $. The determinant of the submatrix in the dashed box above corresponds to  $\det \frac{\partial \phi^2_t}{\partial y}$ in the definition of $P$. Hence, for any initial condition in the horizontal plane $(x,y,0)$, we have that $P= [\pi/2, \infty)$. Geometrically, the dynamics of the horizontal plane $\{(x,y,0); x, y \in \R \}$ is simply a rotation by $f(t)$ around the $z$-axis, with a constant velocity increasing in the $z$-coordinate. Hence, at time $t\geq \frac{\pi}{2}$, the image of $e_1$ ($\in \R^p$, with $p=1$) intersects the plane $\{(0,y,z); y, z \in \R \}= \R^2$, the vertical foliation.}

\end{example}

\bigskip

\section{General case}

\subsection{`Stop and go' technique and preliminaries}

The technique of the previous section (i.e., via Marcus equation) does not allow one to generate a $p$-decomposable flow $\widetilde{\varphi}_t$ of diffeomorphisms close to the original flow $\varphi_t$ when the vector fields do not commute with each other. In fact, the jumps of Marcus equation occur in the direction of the deterministic flows (they do not reach, for example, the directions of their Lie brackets). In other words, one can not control the proximity of $\widetilde{\varphi}_t$ and $\varphi_t$ just controlling how close $Z_t$ is from $U_t$.

The idea of this section is to propose a technique to generate a process $\widetilde{\varphi}_t$, such that, when necessary, it stops at a certain point, then it jumps, at appropriate times, landing exactly at the original flow $\varphi_t$. This tool, that we call `stop and go' technique, allows one to obtain an analogous result of Proposition 2.3 for the general case.

Using the same notation as before, consider the Stratonovich SDE given by:
\begin{equation}\begin{array}{lll}
\label{eq: original}
 dx_t =  X^0 (x_t)\ dt + \displaystyle \sum_{i=1}^{m} X^i(x_t) \circ dB^i_t,
\end{array}
\end{equation}
whose solution is given by the flow of diffeomorphisms $\varphi_t$, but here, the vector fields not necessarily commute with each other. Given an arbitrary  sequence $T$ of stopping times, i.e.:
$$T: \ \ \ \ 0 = \overline{T}_0 \leq T_1 < \overline{T}_1 \leq T_2 < \overline{T}_2 < \ldots$$
we define the `stop and go' equation, as follows:

\begin{equation}\begin{array}{lll}
\label{eq: stop}
dx_t = \displaystyle \sum_{i=0}^{m} X^i(x_t) \ \square^T \ d U^i_t,
\end{array}
\end{equation}
where, again $U_t = (t, B_t^1, \ldots ,  B_t^m)$ and whose solution must be interpreted, in terms of action of diffeomorphisms, as the process:

\[
\widetilde{\varphi}_t = \begin{cases}
\varphi_t, \ \ \text{if} \ \ t \in [\overline{T}_j, T_{j+1}); \\
\varphi_{T_j}, \ \ \text{if} \ \ t \in [T_j, \overline{T}_{j}).
\end{cases}
\]

\


\subsection{Main result}

 Consider the SDE (\ref{eq: original}) and its corresponding stochastic flow $\varphi_t$. Given an initial condition, for each $\omega \in \Omega$, define the set $P(\omega) = \{t \in \R : \det \frac{\partial \varphi_t^2}{\partial y}(\omega) = 0 \}$, which is a random set, different from the previous section.
 
\begin{proposition}
Given $\varepsilon > 0$, and an open random set $A(\omega) \supseteq P(\omega)$ such that $\mu (A(\omega) \setminus P(\omega))< \varepsilon$, there exists a sequence $T$ of stopping times such that:

\begin{itemize}
\item The process $\widetilde{\varphi}_t$, solution of the `stop and go' equation 
    (\ref{eq: stop}), is p-decomposable for all $t \geq 0$;

\item $\widetilde{\varphi}_t (\omega) = \varphi_t (\omega)$ if $t \notin A(\omega)$.
\end{itemize}

\begin{proof}

Since $0 \notin P (\omega)$ for all $\omega\in \Omega$, we can assume (reducing the set $A (\omega)$, if necessary), that $0 \notin A (\omega)$. We define a random partition of $\R$ in an analogous way of the previous results. Take, for each $\omega$ a continuous  function $F_{\omega}: \R \to [0,1]$ with the property that $F_{\omega}^{-1} (1) = A^C(\omega)$ and $F_{\omega}^{-1} (0) = P(\omega)$, and define:

\begin{itemize}
\item `Green zone' $G_{\omega} = F_{\omega}^{-1} (1)$;
\item `Yellow zone' $Y_{\omega} = F_{\omega}^{-1} \big(\frac{1}{2}, 1 \big)$;
\item `Red zone'  $R_{\omega} = F_{\omega}^{-1} \big[0 , \frac{1}{2}\big].$
\end{itemize}

Define by induction a sequence  of stopping times $T$, as follows: $\overline{T}_0 = 0$, $T_i(\omega) = \inf \ \{ t > \overline{T}_{i-1}(\omega): t \in R_{\omega}\}$ and $\overline{T}_i(\omega) = \inf \ \{ t > T_i(\omega): t \in G_{\omega}\}$.

The sequence $T$ was constructed in a such way that $\det \frac{\partial \widetilde{\varphi}_{t}^2}{\partial y}(\omega)$ does not vanish. So $\widetilde{\varphi}_t$ is p-decomposable for all $t \geq 0$. For the second part of the statement, if $t \notin A(\omega)$ then:
$$t \notin \displaystyle \bigcup_{k \in \N} \big(T_k (\omega), \overline{T_k}(\omega)\big),$$ 
and in this case, $\varphi_t (\omega) = \widetilde{\varphi}_t (\omega)$.

\end{proof}
\end{proposition}

\end{document}